\documentclass[final]{siamltex}
 \begin{document}

\title{ Approximation Theory of Matrix Rank Minimization and Its Application to   Quadratic Equations }

\author{YUN-BIN ZHAO\thanks{School of
Mathematics, University of Birmingham, Edgbaston B15 2TT,
Birmingham,  United Kingdom ({\tt y.zhao.2@bham.ac.uk}).}}

\maketitle

\begin{abstract} Matrix rank minimization problems are gaining a
plenty of recent attention in both mathematical and engineering
fields. This class of problems, arising in various and
across-discipline applications, is known to be NP-hard in general.
In this paper, we aim at providing an approximation theory for the
rank minimization problem, and prove  that a  rank minimization
problem can be approximated to any level of accuracy via continuous
optimization (especially, linear and nonlinear semidefinite
programming) problems. One of the main results in this paper shows
that if   the feasible set of the problem has a minimum rank element
with the least F-norm (i.e., Frobenius norm), then the solution of
the approximation problem converges to the  minimum rank solution of
the original  problem as the approximation parameter tends to zero.
The tractability under certain conditions and convex relaxation of
the approximation problem are also discussed. The methodology and
results in this paper provide a new theoretical basis for the
development of some efficient computational methods for solving rank
minimization problems. An immediate application of this theory to
the system of quadratic equations is presented in this paper.  It
turns out that the condition for such a  system without a nonzero
solution can be characterized by a rank minimization problem, and
thus the proposed approximation theory can be used to establish some
sufficient conditions for the system to possess  only zero solution.
\end{abstract}



\section{Introduction}  Throughout the paper,   let $R^n$  be the
$n$-dimensional Euclidean space,    $R^{m\times n} $ be  the
  $m\times n$ real matrix space, and  $S^n$ be the set of real
symmetric matrices. When $X ,Y \in R^{m\times n},$ we use $\langle
X, Y \rangle =\textrm{tr}(X^T Y) $ to denote the inner product of
$X$ and $Y.$   $\|X\|$ and $\|X\|_F$ denote the spectral norm  and
Frobenius norm of $X$, respectively,  and $\|X\|_*$ stands for the
nuclear norm of $X$ (which is the sum of singular values of $X$).
$A\succeq 0 ~(\succ 0)$ means that $A \in S^n $ is positive
semidefinite (positive definite). Given an $X\in R^{m\times n}$ with
rank $r$, we use $\sigma (X)$ to denote the vector $(\sigma_1(X),
..., \sigma_r(X))$ where $\sigma_1(X)\geq  \cdots\geq \sigma_r(X)>0$
are the singular values of $X$.

Let $C\subseteq R^{m\times n}$ be a closed set. Consider the rank
minimization problem:
\begin{equation} \label{rank}    \textrm{Minimize}  ~\left\{\textrm{rank}(X):
  ~~ X\in C \right\},   \end{equation} which has found many applications in system
control \cite{EGG93, EG94, MP97, M98,  HTS99, F02, FHB04}, matrix
completion \cite{CR08, CT09, WYZ10}, machine learning \cite{ABEV06,
MJCD08}, image reconstruction and distance geometry\cite{LLR95,T00,
SY07, RFP07,  D09}, combinatorial and quadratic optimization
\cite{AV09, ZF10},  to name but a few.  The recent work on
compressive sensing (see e.g. \cite{CRT06, C06, D06}) also
stimulates an extensive investigation of this class of problems. In
many applications, $C$ is defined by a linear map ${\cal A}: R^{
m\times n} \to R^p $ . Two typical situations are
\begin{eqnarray} \label{CCC01}  & C= \{X\in
R^{m\times n}: ~ {\cal A} (X)=b \},  & \\
 &  \label{CCC} C=\{ X\in S^n: ~{\cal A} (X)=b, ~X\succeq 0 \}.  & \end{eqnarray}  Unless $C$
has a very special structure, the problem (\ref{rank}) is difficult
to solve due to the discontinuity and nonconvexity of
$\textrm{rank}(X). $ It is NP-hard since it includes the cardinality
minimization as a special case \cite{N95, RFP07}. The existing
algorithms for (\ref{rank}) are largely heuristic-based, such as the
alternating projection \cite{GB00,D09}, alternating LMIs
\cite{SIG98}, and nuclear norm minimization (see e.g. \cite{F02,
FHB04, RFP07, MGC08, TY09, RXH09}. The idea of the nuclear norm
heuristic is to replace the objective of (\ref{rank}) by the nuclear
norm $\|X\|_*,$ and to solve the following convex optimization
problem:
 \begin{equation} \label{nuclear}   \textrm{Minimize} ~  \{ \|X\|_* :   ~X\in C \}.  \end{equation}
Under some conditions, the solution to the nuclear norm heuristic
coincides with the minimum rank solution (see e.g. \cite{F02, RFP07,
RXH09}). This inspires  an extensive and fruitful study on various
algorithms for solving the nuclear norm minimization problem  \cite
{F02, RFP07,
 MGC08, GM09, TY09, CSPW09, AI10}.
  However,  as pointed out in \cite{RFP07,
RXH09}, the nuclear norm heuristic exhibits a phase transition where
for sufficiently small values of the rank the heuristic always
succeeds, but in the complement of the region, it may fail or never
succeed. While  the nuclear norm $\|X\|_*$ is the convex envelop of
rank(X) on the unit ball $\{X: \|X\|\leq 1\}$ (see \cite{F02,
RFP07}), it   may have a drastic deviation from the rank of $X$
since $\textrm{rank}(X)$ is a discontinuous concave function, and
hence
  it is not a high-quality approximation of
 $\textrm{rank}(X). $ As a result, the true relationship between
 (\ref{rank}) and (\ref{nuclear}) are not known in many situations unless some strong
assumptions such as ``restricted isometry"  hold \cite{RFP07}.

In this paper, we develop a new approximation theory for rank
minimization problems. We first provide a continuous approximation
for $\textrm{rank} (X), $ by which $\textrm{rank} (X)$ can be
approximated to any prescribed accuracy, and can be even computed
exactly by
 a suitable choice of the approximation parameter. Based
on this fact, we prove that
 (\ref{rank}) can be approximated to any level
of accuracy by a continuous optimization problem, typically, a
structured linear/nonlinear semidefinite programming (SDP) problem.
One of our main results shows that when the feasible set is of the
form (\ref{CCC}), and if it contains a minimum rank element with the
least F-norm (i.e. Frobenius norm), then the rank minimization
problem can be approximated to any level of accuracy via  an SDP
problem, which is computationally tractable. A key feature of the
proposed approximation approach is that the inter-relationship
between (\ref{rank}) and its approximation counterpart can  be
clearly displayed in many situations. The approximation theory
presented in this paper, aided with modern convex optimization
techniques, provides a theoretical basis for (and can directly lead
to) both new heuristic and  exact algorithms for tackling rank
minimization problems.

To demonstrate an application of the approximation theory, let us
consider the system
\begin{equation}\label{quadratic}  x^T
A_i x=0, ~i=1, ...,m, ~x\in R^n,  \end{equation} where $A_i\in S^n,
i=1,..., m.$  A fundamental question associated with
(\ref{quadratic}) is: when is `$x=0'$ the only solution to
(\ref{quadratic})?  The study of this question (e.g. \cite{F37, D43,
B61, U79, HUT02}) can be dated back to the late 1930s. For $m=2$ and
$n\geq 3,$ the answer to the question is well-known: \emph{$0$ is
the only solution to $x^T A_1 x=0, x^T A_2 x =0 $  if and only if
$\mu_ 1 A_1+\mu_2 A_2 \succ 0 $ for some $\mu_1, \mu_2\in R.$}
However, this result is not valid for $n=2,$ or for $m\geq 3.$ In
fact, the condition
\begin{equation} \label{positive}  \sum_{i=1}^m \mu_i A_i
\succ 0 ~\textrm{ for some } \mu_1, ..., \mu_m\in R
\end{equation} implies that $0$ is the only solution to
(\ref{quadratic}), but the converse is not true in general. When
$n=2$ and/or $m\geq 3$, the sufficient condition (\ref{positive})
may be too strong. Thus finding a mild sufficient condition for the
system (\ref{quadratic}) with only zero solution is posted as an
open problem in \cite{HU07}. We first show that the study of this
problem can be transformed equivalently as a rank minimization
problem, based on which we use the proposed approximation theory,
together with the SDP relaxation and duality theory, to  establish
some general sufficient conditions for the system with only zero
solution.

 This paper is organized as follows. In section 2, an approximation function of $\textrm{rank}(X)$ (and thus
 an approximation model for the rank minimization problem) is introduced,
 and some intrinsic properties of this function are shown.  In section
 3, reformulations and modifications of the approximation
 counterpart of the rank minimization problem
 are discussed, and their proximity to the original
 problem is also proved.  The application of the
approximation theory to the system of quadratic equations has been
demonstrated in section 4. Conclusions  are given in the last
section.

\section{Generic approximation of rank minimization}

The objective of this section is to provide an approximation theory
that can be applied  to general rank minimization problems, without
involving a specific structure of the feasible set which is only
assumed to be a closed set (and bounded when necessary, but not
necessarily convex).  In order to get an efficient approximation of
the problem (\ref{rank}), it is natural to start with a sensible
approximation of $\textrm{rank}(X).$  Let us consider the function
$\phi_\varepsilon: R^{m\times n} \to R$  defined by
\begin{equation} \label{phi} \phi_\varepsilon(X) =  \textrm{tr}\left(X(X^TX+\varepsilon I)^{-1} X^T
\right) , ~~ \varepsilon >0  .
\end{equation}
The first result below claims that the rank of a matrix can be
approximated (in terms of $\phi_\varepsilon$) to any prescribed
accuracy, as long as the parameter $\varepsilon $ is suitably
chosen.

 \textbf{Theorem 2.1.} \emph{ Let $X\in R^{m\times n} $ be a matrix
with $\textrm{rank}(X)=r,$ and $\phi_\varepsilon$ be defined by
(\ref{phi}).  Then for every $\varepsilon>0,$
 \begin{equation} \label{phi+} \phi_\varepsilon(X) =\sum_{i=1}^r \frac{(\sigma_i(X))^2}{
(\sigma_i(X))^2 +\varepsilon},
\end{equation}  where   $\sigma_i(X)$'s are the singular
values of $X, $ and the following relation holds:
\begin{equation} \label{relation} ~~~ 0 \leq \textrm{rank}(X)-  \phi_\varepsilon(X)=    \sum_{i=1}^r
\frac{\varepsilon}{(\sigma_i(X))^2+\varepsilon} \leq   \varepsilon
\sum_{i=1}^r \frac{1}{(\sigma_i(X))^2}  ~~ \textrm{ for all }
\varepsilon >0. \end{equation}}

\emph{Proof. } Let $X=U\Sigma V^T$ be the full singular value
decomposition, where $ U, V$ are orthogonal matrices with dimensions
$m$ and $n$, respectively, and the matrix $\Sigma = \left(
                                                           \begin{array}{cc}
                                                             \textrm{diag} (\sigma (X)) & 0_{r\times (n-r)} \\
                                                             0_{(m-r)\times r} & 0_{(m-r)\times (n-r)}\\
                                                           \end{array}
                                                         \right)
$ where $0_{p\times q}$ denotes the $p\times q $  zero matrix. Let
$\sigma^2(X)$ denote the vector $((\sigma_1(X))^2,
     ..., (\sigma_r(X ))^2).$  Note that
$$ X^TX+\varepsilon I  =  V(\Sigma^T \Sigma)V^T +\varepsilon I =V \left(
                                                                    \begin{array}{cc}
                                                                      \textrm{diag} (\sigma^2 (X))+\varepsilon I_r & 0\\
                                                                      0 & \varepsilon  I_{n-r}\\
                                                                    \end{array}
                                                                  \right)V^T,
                                                                  $$
   where $I$ is partitioned into two small identity matrices $I_r$ and $I_{n-r}.$    Thus, we have
     \begin{eqnarray*}   \phi_\varepsilon (X)  & = &  \textrm{tr}\left(X(X^TX+\varepsilon I)^{-1} X^T
     \right)    \\
     & = & \textrm{tr}\left(U\Sigma \left(\begin{array}{cc} \textrm{diag} (\sigma^2 (X))+\varepsilon I_r & 0\\
0 & \varepsilon  I_{n-r}\\
\end{array}
\right)^{-1} \Sigma^T U^T \right) \\
& = & \textrm{tr}\left( \left(\begin{array}{cc} \textrm{diag} (\sigma^2 (X))+\varepsilon I_r & 0\\
0 & \varepsilon  I_{n-r}\\
\end{array}
\right)^{-1} (\Sigma^T \Sigma)  \right)
 \\
 & = & \textrm{tr}\left(\left(\begin{array}{cc} \textrm{diag} (\sigma^2 (X))+\varepsilon I_r & 0\\
0 & \varepsilon  I_{n-r}\\
\end{array}
\right)^{-1}    \left(\begin{array}{cc} \textrm{diag} (\sigma^2 (X)) & 0 \\
                                                             0  & 0 \\
                                                           \end{array}
                                                         \right)      \right) \\
& = & \sum_{i=1}^r \frac{(\sigma_i(X))^2}{ (\sigma_i(X))^2
+\varepsilon}.
\end{eqnarray*}
Clearly, $ \phi_\varepsilon (X) \leq r=\textrm{rank}(X)  $ for all
$\varepsilon >0.$  Note that
$$ \textrm{rank}(X)-\phi_\varepsilon (X) = \sum_{i=1}^r \left(1- \frac{(\sigma_i(X))^2}{
\sigma_i(X))^2 +\varepsilon}\right)= \sum_{i=1}^r
\frac{\varepsilon}{ (\sigma_i(X))^2 +\varepsilon} \leq  \sum_{i=1}^r
\frac{\varepsilon}{ (\sigma_i(X))^2 }. $$ Thus the inequality
(\ref{relation}) holds.  $ \endproof$

From the above result we have $\phi_\varepsilon (X)\leq
\textrm{rank}(X)$ and $\lim_{\varepsilon \to 0} \phi_\varepsilon (X)
= \textrm{rank} (X). $ So, we immediately have the following
corollary.

\textbf{Corollary 2.2.} \emph{For every matrix $X\in R^{m\times n}$,
there exists accordingly a number
 $\varepsilon^* >0$ such that $ \textrm{ rank}({X})   = \lceil
 \phi_\varepsilon (X)\rceil $ for all $ \varepsilon \in (0,
 \varepsilon^*].  $
}

 This suggests  the following scheme which
requires only a finite number of iterations to find the exact rank
of $X$: \emph{Step 1}. Choose a small number $\varepsilon>0; $
\emph{Step 2}.  Evaluate  $\phi_\varepsilon (X) $ at $X; $
\emph{Step 3}.  Round up the value of $\phi_\varepsilon (X)$  to the
nearest integer; \emph{Step 4}. Set $\varepsilon \leftarrow \beta
\varepsilon $ where $\beta\in (0,1) $ is a given constant, and
repeat the steps 2-4 above.

The threshold $\varepsilon^*$ in Corollary 2.2 depends on  $X.$ This
can be seen clearly from the right-hand side of (\ref{relation}).
However, the next theorem shows that over the optimal solution set
of (\ref{rank}) the approximation is uniformed. Before stating this
result, we first show that the optimal solution set   of
(\ref{rank}) is closed. Note that, in general, the set $\{X\in C:
\textrm{rank}(X) =r\} $ is not closed.

\textbf{Lemma 2.3.} \emph{Let $C$ be a closed set in $R^{m\times n}.$ Then
 the level set $ \{X\in C: \textrm{rank}(X)
\leq r\}$ is closed for any given number $r\geq 0. $ In particular,
 the optimal solution set  of (\ref{rank}), i.e.,
$ C^*= \{X\in C: ~ \textrm{rank}(X)= r^*\}   $ is closed, where
$r^*(=\min\{\textrm{rank}(X): X\in C\}) $ is the minimum rank.}

\emph{Proof.}  Suppose that $\{X^k\}\subseteq \{X\in C:
\textrm{rank}(X) \leq r\}$  is a sequence convergent to $X^0$ in the
sense that $ \|X^k-X^0\|\to 0 $ as $k\to \infty.$ Let $r^0
=\textrm{rank}(X^0) $ and $\sigma_1(X^0) \geq \cdots \geq
\sigma_{r^0}(X^0)
>0 $ be the nonzero singular values of $X^0$. Note that the singular
value is continuously dependent on the entries of the matrix. It
implies that for sufficiently large $k$, $X^k$ has at least $r^0$
nonzero singular values. Thus $ \textrm{rank}(X^0) \leq
\textrm{rank}(X^k) \leq r
 $ for all sufficiently large $k. $ This together with the closedness of $C$  implies that
  $X^0 \in \{X\in C: \textrm{rank}(X)\leq r\},$
and thus  the level set of $ \textrm{rank}(X)$  is closed.
Particularly, it implies that  the optimal solution set $\{X\in C:~
\textrm{rank}(X) = r^*\}= \{X\in C:~ \textrm{rank}(X) \leq r^*\}$ is
closed.  $  \endproof$

We now show that the function $\textrm{rank}(X) $ can be uniformly
approximated by $\phi_\varepsilon (X) $ over the optimal solution
set of (\ref{rank}), in the sense that the right-hand side of
(\ref{relation}) is independent of the choice of $X^*.$

 \textbf{Theorem 2.4.} \emph{If the optimal solution set, denoted by $C^*$, of (\ref{rank}) is
 bounded, then there exists a constant $\delta >0$ such that for any
 given $\varepsilon>0$ the inequality
 $$ \phi_\varepsilon (X^*) \leq  rank(X^*) \leq \phi_\varepsilon
 (X^*)+\varepsilon \left(\frac{\min\{m,n\}}{\delta^2}\right) $$  holds for all $ X^*\in
 C^*.$}

\emph{Proof.} Let $r^*$ be the minimum rank of (\ref{rank}). Then
$r^*= \textrm{rank}(X^*)$ for all $X^*\in C^*.$ Let $ \sigma_{r^*}
(X^*)$  denote the smallest nonzero singular value of $X^*, $ and
denote $$\sigma_{min} = \min\{\sigma_{r^*} (X^*): ~X^*\in C^*\}.$$
We now prove that $\sigma_{min} >0. $ Indeed, if $\sigma_{min}=0$,
then there exists a sequence $\{X^{*}_k\} \subseteq C^*$ such that
$\sigma_{r^*} (X^{*}_k) \to 0.$ Since $C^*$ is bounded, passing to a
subsequence if necessary we may assume that $X^{*}_k\to \widehat{X}.
$ Thus, $\sigma_{r^*} (\widehat{X}) =0,$ which implies that $
\textrm{rank}(\widehat{X})< r^*$, contradicting to the closedness of
$C^*$ (see Lemma 2.3).  Therefore, we have $\sigma_{min}>0.$ Let
$\delta >0$ be a constant satisfying $\delta\leq \sigma_{min}.$  By
(\ref{relation}), we have
$$ \textrm{rank}(X^*) -\phi_\varepsilon (X^*) \leq \varepsilon \sum_{i=1}^{r^*} \frac{1}{(\sigma_i(X^*))^2} \leq
 \varepsilon  \frac{r^* }{(\sigma_{r^*}(X^*))^2}  \leq   \varepsilon \left(\frac{\min\{m,n\}}{\delta^2}\right),  $$
as desired.  ~ $  \endproof$

It is easy to see from (\ref{phi}) that   $\phi_\varepsilon (X)$ is
continuous with respect to $(X,\varepsilon) $ over the set
$R^{m\times n} \times (0, \infty).$ From Theorem 2.1 and Corollary
2.2, we see that   the problem (\ref{rank}) can be approximated by a
continuous optimization problem with  $\phi_\varepsilon.$
 In fact, by replacing $\textrm{rank}(X)$ by $\phi_\varepsilon (X), $
 we obtain the following approximation problem of (\ref{rank}):
  \begin{equation} \label{phi-mini} \begin{array}{cl} \textrm{Minimize}  & \phi_\varepsilon (X)= \textrm{tr}\left(X(X^TX+\varepsilon I)^{-1}X^T\right)  \\
\textrm{s.t.} & X\in C
\end{array}
\end{equation}
where $\varepsilon>0$ is a given parameter. From an approximation
point of view,  some natural questions arise:  Does the optimal
value (solution) of (\ref{phi-mini}) converges to a minimum rank
(solution) of (\ref{rank}) as $\varepsilon\to 0$? How can we solve
the problem (\ref{phi-mini}) efficiently, and when this  problem is
computationally tractable? The remainder of this section and the
next section  are devoted to answering these questions.

For the convenience of the later analysis, we use notation $\phi_0
(X)= \textrm{rank}(X). $  Before we prove the main result of this
section, let us first prove the semicontinuity  of the function
$\phi_\varepsilon (X) $ at the boundary point $\varepsilon=0$,

\textbf{Lemma 2.5.} \emph{ With respect to $(X, \varepsilon), $ the
function $\phi_\varepsilon (X)$ is continuous everywhere in the
region $ R^{m\times n} \times (0, \infty), $  and it is lower
semicontinuous at $(X, 0),$ i.e.,
$$\liminf_{(Y, \varepsilon)\to (X, 0)} \phi_\varepsilon (Y) \geq \phi_0(X) =\textrm{rank}(X).
$$ }

\emph{Proof.} The continuity of $\phi_\varepsilon$ in $ R^{m\times
n} \times (0, \infty) $ is obvious. We only need to prove  its lower
semicontinuity at $(X, 0). $  Let $ \widehat{X}$ be an arbitrary
matrix in $ R^{m\times n} $ with $\textrm{rank}(\widehat{X}) =r.$
Suppose that $X  \to \widehat{X}. $ Then it is easy to see that
\begin{equation} \label{SSS} (\sigma_1(X), ..., \sigma_r(X)) \to (\sigma_1(\widehat{X}),
..., \sigma_r(\widehat{X}))
>0,\end{equation}  and $ \sigma_i(X) \to 0$ for $i\geq r+1.$ This implies that
$\textrm{rank}(X) \geq \textrm{rank}(\widehat{X}) $ as long as $X$
is sufficiently close to $\widehat{X}.$  By (\ref{phi+}), we have
\begin{eqnarray} \label{semicont} \phi_\varepsilon (X) -
\phi_{0}(\widehat{X})  & = &  \phi_\varepsilon (X) -
\textrm{rank}(\widehat{X}) \nonumber
\\ &  = &  \sum_{i\geq r+1}
\frac{(\sigma_i (X))^2}{(\sigma_i (X))^2+\varepsilon } +\sum_{i=1}^r
\left(\frac{(\sigma_i (X))^2}{(\sigma_i (X))^2+\varepsilon} - 1
\right).
 \end{eqnarray}
 It is not difficult to see that  when $(X, \varepsilon)\to (\widehat{X}, 0)$, the right-hand side
 of  (\ref{semicont})  does not necessarily tend to zero, when $(\sigma_i (X))^2$ in the first
 term of the right-hand side of (\ref{semicont}) tends to zero no faster than that of $\varepsilon.$   For instance, let
 $(\sigma_1 (\widehat{X}), ..., \sigma_r (\widehat{X}))  =
 (1,...,1)$, and consider the sequence $X^k \to \widehat{X} $ where $X^k$ satisfies that
  $\textrm{rank}(X^k) = p >r ,$  $ (\sigma_1 (X^k), ..., \sigma_r (X^k))  =
 (1,...,1)$ and $ (\sigma_{r+1} (X^k), ..., \sigma_p (X^k))  =
 (1/k,...,1/k). $  Setting  $\varepsilon_k
 =\frac{1}{k^2}  $  and substituting $(X^k, \varepsilon^k)$ into  (\ref{semicont})  yields
$$ \phi_{\varepsilon_k} (X^k) - \phi_{0}(\widehat{X})
 =  \frac{1}{2} (p-r) - \frac{r}{1+k^2}   \to   \frac{1}{2} (p-r) >0 ~ \textrm{ as } k\to \infty.$$\\
So, $\phi_\varepsilon(X)$ is not necessarily continuous at
$\varepsilon =0.$ However, from (\ref{semicont}) we see that
$$ \phi_\varepsilon (X) - \phi_{0}(\widehat{X}) \geq   \sum_{i=1}^r
\left(\frac{(\sigma_i (X))^2}{(\sigma_i (X))^2+\varepsilon} - 1
\right), $$  where $r = \textrm{rank} (\widehat{X}).$ By
(\ref{SSS}), the right-hand side of the above goes to zero as $X\to
\widehat{X} $ and $ \varepsilon \to 0.$ It follows that
$$ \liminf_{(X, \varepsilon)\to (\widehat{X}, 0)} \phi_\varepsilon (X)
\geq \liminf_{(X, \varepsilon)\to (\widehat{X}, 0)}
\left(\phi_0(\widehat{X})+ \sum_{i=1}^r \left(\frac{(\sigma_i
(X))^2}{(\sigma_i (X))^2+\varepsilon} - 1 \right) \right)=
\phi_0(\widehat{X}).$$ The proof is complete. ~ $  \endproof$

We now prove the main result of this section, which shows that  the
rank minimization over a bounded feasible set can be approximated
with (\ref{phi-mini}) to any level of accuracy.

\textbf{Theorem 2.6.} \emph{Let $C$ be a closed set in $R^{m\times
n}.$ Let $r^*$ be the  minimum rank  of (\ref{rank}) and for given
$\varepsilon>0$, let $\phi^*_\varepsilon $ and $X(\varepsilon)$ be
the optimal value  and  an optimal solution of (\ref{phi-mini}),
respectively. Then
\begin{equation} \label{underestimate}
\phi^*_\varepsilon \leq r^* ~ \textrm{ for any } \varepsilon >0.
\end{equation}
 Moveover, when $C$ is bounded, then
\begin{equation} \label{limit} \lim_{\varepsilon \to 0} \phi^*_\varepsilon  =r^*, \end{equation}  and  any
accumulation point of $X(\varepsilon)$, as $\varepsilon\to 0,$ is a
minimum rank solution of (\ref{rank}).}

\emph{Proof.} Since $X(\varepsilon)$ is an optimal solution to
(\ref{phi-mini}), we have
$$\phi_\varepsilon ^* =\phi_\varepsilon
(X(\varepsilon)) \leq \phi_{\varepsilon} (X) ~~ \textrm{ for all }
X\in C.$$ Particularly, any optimal solution of (\ref{rank})
satisfies the above inequality. So, $\phi^*_\varepsilon  \leq
\phi_\varepsilon (X^*) \leq \textrm{rank}(X^*)= r^* $ where the
second inequality follows from (\ref{relation}), and $\varepsilon$
can be any positive number. Thus (\ref{underestimate}) holds,  and
\begin{equation} \label{limsup} \limsup_{\varepsilon\to 0} \phi^*_\varepsilon \leq r^*. \end{equation}  On
the other hand, since  $\phi^*_\varepsilon \geq 0,$ the number
$\underline{r} =\liminf_{\varepsilon\to 0} \phi^*_\varepsilon$ is
finite. Without loss of generality, assume that the sequence
$\{\phi^*_{\varepsilon_k}\} $, where $\varepsilon_k\to 0$ as $k\to
\infty$,  converges to $\underline{r}. $ Note that
$X(\varepsilon_k)$ is a minimizer of (\ref{phi-mini}) with
$\varepsilon=\varepsilon_k,$ i.e., $\phi^*_{\varepsilon_k} =
\phi_{\varepsilon_k} (X(\varepsilon_k)).$ When $C$ is bounded, the
sequence $\{X(\varepsilon_k)\}$ is bounded. Passing to a subsequence
if necessary, we assume that $ X(\varepsilon_k) \to X_0$ as $k\to
\infty.$  Clearly, $X_0\in C$ since $C$ is closed, and hence
$\textrm{rank}(X_0) \geq r^*.$ Therefore,
$$ \underline{r} =\lim_{k\to \infty} \phi^*_{\varepsilon_k}
=\lim_{k\to \infty} \phi_{\varepsilon_k} (X_{\varepsilon_k}) \geq
\liminf_{(X,\varepsilon)\to (X_0, 0)} \phi_\varepsilon (X) \geq
\phi_0 (X_0),$$ where the last inequality follows from Lemma 2.5.
Thus, $ \underline{r} \geq \phi_0 (X_0)=\textrm{rank}(X_0)\geq r^*,
$ which together with (\ref{limsup}) implies  (\ref{limit}).

We now prove that any accumulation point of $X(\varepsilon)$ is a
minimum rank solution of (\ref{rank}). Let $\widehat{X}$ (with
$\textrm{rank}(\widehat{X}) =\widehat{r} $) be an arbitrary
accumulation point of $X(\varepsilon)$, as $\varepsilon \to 0.$ We
now prove that $\widehat{X}$ is a minimum rank solution to
(\ref{rank}), i.e., $ \widehat{r} =r^* .$ Consider a convergent
sequence
  $X(\varepsilon_k) \to \widehat{X} $ where $ \varepsilon_k \to 0. $
  Then by
 (\ref{underestimate}) and (\ref{phi+}), we have
\begin{eqnarray*} r^* \geq \phi^*_{\varepsilon_k}= \phi_{\varepsilon_k}
(X(\varepsilon_k))
 & = & \sum_{i=1}^{\widehat{r}} \frac{(\sigma_i (X(\varepsilon_k)))^2 }{ (\sigma_i (X(\varepsilon_k)))^2
 +\varepsilon_k} + \sum_{i> \widehat{r}} \frac{(\sigma_i (X(\varepsilon_k)))^2 }{ (\sigma_i (X(\varepsilon_k)))^2
  +\varepsilon_k} \\
   & \geq  & \sum_{i=1}^{\widehat{r}} \frac{(\sigma_i (X(\varepsilon_k)))^2 }{ (\sigma_i (X(\varepsilon_k)))^2
 +\varepsilon_k} \to \sum_{i=1}^{\widehat{r}} \frac{(\sigma_i (\widehat{X}))^2 }{ (\sigma_i (\widehat{X}))^2
 +0} =\textrm{rank} (\widehat{X}).
  \end{eqnarray*}
Thus, any accumulation point of $X(\varepsilon)$ is a minimum rank
solution to (\ref{rank}). ~  $  \endproof$

Since $r^*$ is integer, by (\ref{underestimate}) and (\ref{limit}),
we immediately have the corollary below.

\textbf{Corollary 2.7.} \emph{ Let $r^*$ and $\phi_\varepsilon^* $
be defined as in Theorem 2.6. If $C \subset R^{m\times n} $ is
bounded and closed, then there exists a number $\delta >0$ such that
$r^*= \left \lceil  \phi_\varepsilon^*  \right  \rceil $ for all $
\varepsilon \in (0, \delta]. $}

 The results above provide a theoretical basis for
developing new approximation algorithms for rank minimization
problems. Such an algorithm can be a heuristic method for general
rank minimization, and can be  an exact method as indicated by
Corollary 2.7. From Theorem 2.6, the set $\{X(\varepsilon)\}$ can be
viewed as a trajectory leading to the minimum rank solution set of
(\ref{rank}), and thus it is possible to construct a continuation
type method (e.g. a path-following method) for rank minimization
problems. In the next section, we are going to discuss how and when
the approximation problem (\ref{phi-mini}) can be efficiently dealt
with from the viewpoint of computation. We prove that under some
conditions problem (\ref{phi-mini}) can be either reformulated or
relaxed as a tractable optimization problem, typically an SDP
problem.

\section{Reformulation of the approximation problem (\ref{phi-mini})}
The main result in last section shows that  if $\varepsilon $ is
small enough, the optimal value of the rank minimization problem can
be obtained precisely by solving (\ref{phi-mini}) just once, and the
solution $X(\varepsilon) $  of (\ref{phi-mini}) is an approximation
to the optimal solution of (\ref{rank}). If the problem
(\ref{phi-mini}) with a prescribed $\varepsilon>0 $ fails to
generate the minimum value of (\ref{rank}), we can reduce the value
of $\varepsilon$ and solve (\ref{phi-mini}) again. By Corollary 2.7,
the minimum rank of (\ref{rank}) can be obtained by solving
 (\ref{phi-mini}) up to a finite number of times. Thus, roughly speaking, solving a rank minimization problem amounts
to solving a continuous optimization problem defined by
(\ref{phi-mini}).

In this section, we concentrate on the problem (\ref{phi-mini}) to
find out when and how it can be solved efficiently. To this end, we
investigate its equivalent formulations together with some useful
variants. By doing so, we take into account the structure of $C$
when necessary.  Let us start with the reformulation of
(\ref{phi-mini}).

Introducing a variable $Y\in S^m,$  we first note that
(\ref{phi-mini}) can be written as the following
 nonlinear semidefinite programming problem:
\begin{equation} \label{phi-mini-A}  \textrm{Minimize} \left \{
\textrm{tr}(Y): ~ Y \succeq X(X^TX+\varepsilon I)^{-1} X^T ,
  ~ X\in C\right\}.
\end{equation}
It is easy to see that if $(Y^*, X^*)$ is an optimal solution to
(\ref{phi-mini-A}), then   \begin{equation}\label{YYY} Y^*=
X^*((X^*)^TX^*+\varepsilon I)^{-1} (X^*)^T. \end{equation}   Thus,
we conclude that $X^*$ is an optimal solution to (\ref{phi-mini}) if
and only if $(Y^*, X^*)$ is an optimal solution to
(\ref{phi-mini-A}) where $Y^*$ is given by (\ref{YYY}).  By Schur
complement theorem, the problem (\ref{phi-mini-A}) can be further
written as
\begin{equation} \label{phi-mini-B}   \textrm{Minimize}
\left\{\textrm{tr}(Y):  ~~  \left(
           \begin{array}{cc}
             Y & X \\
             X^T &  X^TX+ \varepsilon I \\
           \end{array}
         \right)
 \succeq  0, ~~ X\in C\right\},
\end{equation}
which remains a  nonlinear SDP problem. We now introduce the
variable $Z= X^TX, $ which implies that $Z$ is the optimal solution
to the problem $ \min_{Z} \{\textrm{tr}(Z): Z\succeq X^T X\}. $ By
Schur complement theorem again, $ Z\succeq X^T X $ is nothing but $
\left(
\begin{array}{cc}
                    I & X  \\
                    X^T & Z \\
                  \end{array}
                \right) \succeq 0. $
                So the problem (\ref{phi-mini}) can be written
                 exactly as a
 bilevel SDP problem:
\begin{eqnarray} \label{bilevel}   & \min_{ Y\in S^m,\widetilde{Z}\in S^n, X\in C}   &  \textrm{tr}(Y)  \nonumber \\
 & \textrm{s.t.} &   \left(
           \begin{array}{cc}
             Y & X \\
             X^T &  \widetilde{Z}+ \varepsilon I \\
           \end{array}
         \right)
 \succeq  0,  \\
&  &   \widetilde{Z}= \textrm{arg}\min_{Z\in S^n} \left\{
\textrm{tr}(Z): \left(
                  \begin{array}{cc}
                    I & X  \\
                    X^T & Z \\
                  \end{array}
                \right) \succeq 0 \right \}.   \nonumber
\end{eqnarray}
From the discussion above, we conclude that (\ref{phi-mini}) is
equivalent to the nonlinear SDP problem (\ref{phi-mini-B}), and is
equivalent to the \emph{linear bilevel SDP problem} (\ref{bilevel}).
As a result, by Theorem 2.6,  \emph{the rank minimization over a
bounded feasible set is equivalent to the linear bilevel SDP problem
of the form (\ref{bilevel}).} Thus, the level of difficulty for rank
minimization can be understood from the perspective of its linear
bilevel SDP counterpart. It is worth mentioning that the bilevel
programming (in vector form) has been long studied (e.g.
\cite{LPR96}), but to our knowledge the bilevel SDP problem remains
a new topic  so far. The analysis above shows that a bilevel SDP
model does arise from rank minimization.  However, both
(\ref{phi-mini-B}) and (\ref{bilevel}) are not convex problems, and
hence they are not computationally tractable in general.

This motivates us to consider the next approximation model which can
be viewed as a variant of (\ref{phi-mini}). The difficulty of
(\ref{phi-mini-B}) and (\ref{bilevel}) lies in the hard equality
$Z=X^T X. $  An immediate idea is to relax it to  $Z \succeq X^T X,
$ yielding the
  problem:
$$   \textrm{Minimize} \left\{\textrm{tr}(Y):
~~   \left(
           \begin{array}{cc}
             Y & X \\
             X^T &  Z+ \varepsilon I \\
           \end{array}
         \right)
 \succeq  0,    \left(
                  \begin{array}{cc}
                    I & X  \\
                    X^T & Z \\
                  \end{array}
                \right) \succeq 0 ,
 ~~ X\in C
 \right\}
$$
which is a convex problem if $C$ is convex, and an SDP problem if
$C$ is defined by (\ref{CCC01}) or (\ref{CCC}). However, for any
given $X\in C$ and any number $\beta>0$, the point $( Y=\beta I, Z=
\alpha I) $ is feasible to the above problem when  $\alpha
>0 $ is sufficiently large. So the
optimal value of the above problem is always zero, providing nothing
 about the minimum rank of the original problem
(\ref{rank}). This happens since $Z$ gains too much freedom while
$Z= X^T X$ is relaxed to $Z \succeq  X^TX. $  Thus the value $
\textrm{tr}(Y) = \textrm{tr}(X (Z+\varepsilon I)^{-1} X^T)$  may
significantly deviate from $\phi_\varepsilon (X)
(\approx\textrm{rank} (X)). $ To avoid this, some driving force
should be imposed on $Z$ so that it is near (or equal) to $X^T X. $

Motivated by this observation, we consider the following problem in
which  a `penalty' term is introduced into  the objective:
 \begin{equation} \label{phi-mini-C} \begin{array}{cl} \textrm{Minimize}  & \textrm{tr}(Y) + \frac{1}{\gamma} \textrm{tr}(Z)\\
s.t. &   \left(
           \begin{array}{cc}
             Y & X \\
             X^T &  Z+ \varepsilon I \\
           \end{array}
         \right)
 \succeq  0,    \left(
                  \begin{array}{cc}
                    I & X  \\
                    X^T & Z \\
                  \end{array}
                \right) \succeq 0 ,
 ~~ X\in C,
\end{array}
\end{equation}
where $\gamma$ is a positive number. The term $\frac{1}{\gamma}
\textrm{tr}(Z)$ acts as a penalty when $Z (\succeq X^T X) $ is
deviated away from $X^T X. $ Since $\textrm{tr}(Z)\geq
\textrm{tr}(X^TX) = \|X\|_F^2,$ this term also drives $ \|X\|_F $ to
be minimized.  Note that when $Z$ is driven  near to $ X^TX,$ it is
the first term $\textrm{tr}(Y)$ of the objective that approximates
the rank of $X,$ and returns the approximate value of
$\textrm{rank}(X).$ The advantage of the approximation model
(\ref{phi-mini-C}) is that it is an SDP problem when $C$ is defined
by linear constraints (such as (\ref{CCC01}) or (\ref{CCC})), and
hence it is computationally tractable. In what follows, we
concentrate on this model and prove that under some conditions the
rank minimization can be approximated by (\ref{phi-mini-C}) to any
level of accuracy.

\textbf{Theorem 3.1.} \emph{Let $C  $ be a bounded, closed set in
$R^{m\times n}.$ Suppose that $C$ contains a minimum rank element
$X^*$ with the least F-norm, i.e., $ \textrm{rank}(X^*)\leq
\textrm{rank}(X),
 ~ \|X^*\|_F\leq \|X\|_F  \textrm{  for all  }X\in C. $
  Let $(Y_{\varepsilon,
\gamma}, Z_{\varepsilon, \gamma}, X_{\varepsilon, \gamma})$ denote
the optimal solution  of the   problem (\ref{phi-mini-C}).   Then $
\textrm{tr}(Y_{\varepsilon, \gamma}) \leq \phi_\varepsilon (X^*)
\leq rank(X^*) $ for all $(\varepsilon, \gamma)>0$ and
$$ \lim_{(\varepsilon, \gamma)\to 0, \frac{\gamma}{\varepsilon} \to
0} tr(Y_{\varepsilon, \gamma}) = r^* = \textrm{rank}(X^*) , ~
\lim_{(\varepsilon, \gamma)\to 0, \frac{\gamma}{\varepsilon} \to 0}
tr(Z_{\varepsilon, \gamma})= \|X^*\|_F^2,  $$
 and the
sequence $\{X_{\varepsilon, \gamma} \} $  converges to the set of
   minimum rank solutions of (\ref{rank}), as $(\varepsilon,
\gamma)\to 0$ and $ \frac{\gamma}{\varepsilon} \to 0. $ }

\emph{Proof.} Since $C$ is bounded and closed, the sequence
$\{X_{\varepsilon, \gamma}\} $ has at least one accumulation point,
and any such an accumulation point is in $C.$ Let $X^0$ be an
arbitrary accumulation point of the sequence $\{X_{\varepsilon,
\gamma}\} $ as $(\varepsilon, \gamma)\to 0$ and $
\frac{\gamma}{\varepsilon} \to 0.$  Without loss of generality, we
assume that $X_{\varepsilon, \gamma} \to X^0,  $ as $(\varepsilon,
\gamma)\to 0$ and $ \frac{\gamma}{\varepsilon} \to 0.$  By Schur
complement and the structure of the problem (\ref{phi-mini-C}), it
is easy to see that for any given $\varepsilon, \gamma
>0$ the optimal solution  $(Y_{\varepsilon,
\gamma}, Z_{\varepsilon, \gamma}, X_{\varepsilon, \gamma})$ of
(\ref{phi-mini-C}) satisfies the following relation
\begin{equation}\label{P1} Y_{\varepsilon, \gamma} = X_{\varepsilon,
\gamma}(Z_{\varepsilon, \gamma}+\varepsilon I)^{-1} X_{\varepsilon,
\gamma}^T, ~~Z_{\varepsilon, \gamma} \succeq X_{\varepsilon,
\gamma}^T X_{\varepsilon, \gamma}.   \end{equation} Let $X^*$ be an
arbitrary minimum rank solution of (\ref{rank}) with the least
F-norm. Then
 the point $(Y^*_\varepsilon, Z^*, X^*),$ where
 $Y^*_\varepsilon= X^*((X^*)^TX^*+\varepsilon)^{-1} (X^*)^T $ and  $Z^*= (X^*)^TX^*,  $ is feasible
  to the problem (\ref{phi-mini-C}). By optimality, we have
\begin{equation}\label{OPTIM}
\textrm{tr}(Y_{\varepsilon, \gamma}) +\frac{1}{\gamma}
\textrm{tr}(Z_{\varepsilon, \gamma}) \leq
\textrm{tr}(Y^*_\varepsilon)+ \frac{1}{\gamma} \textrm{tr}(Z^*) =
\phi_\varepsilon (X^*)+ \frac{1}{\gamma}\|X^*\|_F^2.
\end{equation}
It follows from (\ref{P1}) that \begin{equation} \label{P2}
\textrm{tr}(Z_{\varepsilon, \gamma}) \geq \textrm{tr} (
X_{\varepsilon, \gamma}^T X_{\varepsilon, \gamma} )=
\|X_{\varepsilon, \gamma}\|_F^2 \geq \|X^*\|_F^2~ \textrm{ for all }
\varepsilon, \gamma>0 .\end{equation} Combining (\ref{OPTIM}) and
(\ref{P2}) yields
\begin{eqnarray}\label{P3}  & \textrm{tr}(Y_{\varepsilon, \gamma}) \leq
\phi_\varepsilon (X^*), & \\
& ~~~~~~ \label{P4} 0\leq \textrm{tr}(Z_{\varepsilon, \gamma}) -
\|X^*\|_F^2 \leq \gamma (\phi_\varepsilon (X^*) - \textrm{tr} (
Y_{\varepsilon, \gamma}) ) \leq \gamma \phi_\varepsilon (X^*) \leq
\gamma \min\{m,n\},  & \end{eqnarray}  for all  $(\varepsilon,
\gamma)>0.$ The last inequality of (\ref{P4}) follows from
$\phi_\varepsilon (X^*) \leq \textrm{rank}(X^*) \leq \min\{m,n\}.$
Let $
 X_{\varepsilon, \gamma} =X^0+\Delta_{\varepsilon, \gamma}$ where
  $ \Delta_{\varepsilon, \gamma}\to 0 $ since $X_{\varepsilon, \gamma} \to X^0. $ Then
 $$X_{\varepsilon, \gamma}^TX_{\varepsilon, \gamma} = (X^0)^T X^0 + (\Delta_{\varepsilon, \gamma}^T X^0 +(X^0)^T\Delta_{\varepsilon,
  \gamma} + \Delta_{\varepsilon, \gamma}^T
 \Delta_{\varepsilon, \gamma}) = (X^0)^T X^0+ G (\Delta_{\varepsilon,
 \gamma}),$$
 where $G( \Delta_{\varepsilon,
 \gamma}) = \Delta_{\varepsilon, \gamma}^T X^0 +(X^0)^T\Delta_{\varepsilon,
  \gamma} + \Delta_{\varepsilon, \gamma}^T
 \Delta_{\varepsilon, \gamma}. $
 Thus by (\ref{P1}) we have \begin{equation} \label {P5}  Z_{\varepsilon, \gamma} \succeq X_{\varepsilon,
\gamma}^T X_{\varepsilon, \gamma} = (X^0)^T X^0+ G
(\Delta_{\varepsilon, \gamma}).  \end{equation}  Note that
$\textrm{tr}(G(\Delta_{\varepsilon, \gamma}))  \to 0$ as $
\Delta_{\varepsilon, \gamma} \to 0.$  By (\ref{P4}) and (\ref{P5}),
we have
$$\|X^*\|_F^2=\lim_{(\varepsilon, \gamma)\to 0, \gamma/\varepsilon \to 0} \textrm{tr}(Z_{\varepsilon,
\gamma}) \geq  \lim_{(\varepsilon, \gamma)\to 0, \gamma/\varepsilon
\to 0}  \textrm{tr} ((X^0)^T X^0+ G (\Delta_{\varepsilon, \gamma}))=
\|X^0\|_F^2.$$ Thus, $X^0$ is a least F-norm element in $C.$ On the
other hand, from (\ref{P5}), we see that $
\widehat{\Delta}_{\varepsilon, \gamma} :=  Z_{\varepsilon, \gamma} -
X_{\varepsilon, \gamma}^T X_{\varepsilon, \gamma} \succeq 0. $ Thus,
by (\ref{P5}) and (\ref{P4}) again, we have
$$\|\widehat{\Delta}_{\varepsilon, \gamma}\|\leq  \textrm{tr}(\widehat{\Delta}_{\varepsilon,
\gamma})= \textrm{tr}(Z_{\varepsilon, \gamma})-\|X_{\varepsilon,
\gamma}\|_F^2 \leq \textrm{tr}(Z_{\varepsilon, \gamma}) -
\|X^*\|_F^2 \leq \gamma \min\{m, n\}.
$$
The first inequality above follows from the fact $
\widehat{\Delta}_{\varepsilon, \gamma}  \succeq 0,$ and the  second
follows from  $\|X_{\varepsilon, \gamma}\|_F \geq \|X^*\|_F. $
Therefore,
\begin{equation} \label{limlim} \|\widehat{\Delta}_{\varepsilon, \gamma}\|
/\varepsilon  \to 0, \textrm{ as }(\varepsilon, \gamma)\to 0
\textrm{ and }\gamma/\varepsilon \to 0.
\end{equation}
 When $M\in R^{n\times n}$ and
  $\|M\| <1$, it is well-known that
 $(I+M)^{-1}= I-M+M^2-M^3+ \cdots = I+ \sum_{i=1}^\infty (-1)^{i}
 M^{i} .   $ Thus, for any  $ U,V\in R^{n\times n}$ where $U$ is
 nonsingular,
 if $\|VU^{-1}\|<1$ we have
 \begin{equation}\label{formula} ~~~~~~(U+V)^{-1}= U^{-1} (I+VU^{-1})^{-1} =  U^{-1} + U^{-1}
 \left(\sum_{i=1}^\infty (-1)^{i}
 (VU^{-1})^{i} \right).   \end{equation}
As $(\varepsilon, \gamma)\to 0$ and $ \gamma/\varepsilon \to 0,$  it
follows from (\ref{limlim}) that
$$ \|\widehat{\Delta}_{\varepsilon, \gamma} (X_{\varepsilon,
\gamma}^T X_{\varepsilon, \gamma} +\varepsilon I)^{-1} \| \leq \|
\widehat{\Delta}_{\varepsilon, \gamma}\| \|  (X_{\varepsilon,
\gamma}^T X_{\varepsilon, \gamma} +\varepsilon I)^{-1} \| \leq \|
\widehat{\Delta}_{\varepsilon, \gamma}\| /\varepsilon \to 0. $$
 Thus, substituting $U= X_{\varepsilon, \gamma}^T X_{\varepsilon, \gamma} +\varepsilon
I$ and $V=\widehat{\Delta}_{\varepsilon, \gamma}$ into
(\ref{formula}) yields
\begin{eqnarray*} & &  (X_{\varepsilon, \gamma}^T X_{\varepsilon,
\gamma} +\widehat{\Delta}_{\varepsilon, \gamma}+\varepsilon I)^{-1}
- (X_{\varepsilon, \gamma}^T X_{\varepsilon, \gamma} +\varepsilon
I)^{-1}  \\
&  & =  (X_{\varepsilon, \gamma}^T X_{\varepsilon, \gamma}
+\varepsilon I)^{-1}  \left( \sum_{ i=1}^\infty (-1)^{i}
 \left(\widehat{\Delta}_{\varepsilon, \gamma}
(X_{\varepsilon, \gamma}^T X_{\varepsilon, \gamma} +\varepsilon
I)^{-1} \right)^{i}\right).
\end{eqnarray*}
Note that $ \|(X_{\varepsilon, \gamma}^T X_{\varepsilon, \gamma}
+\varepsilon I)^{-1}\| \leq 1/\varepsilon $  and  $\left\|
\left(I-\varepsilon (X_{\varepsilon, \gamma}^T X_{\varepsilon,
\gamma} +\varepsilon I)^{-1} \right) \right\| \leq 1.$ When $\left\|
 \widehat{\Delta}_{\varepsilon, \gamma} \right\| /\varepsilon <1, $  we have
\begin{eqnarray*}
& & \left\|X_{\varepsilon, \gamma}^T X_{\varepsilon, \gamma}
\left((X_{\varepsilon, \gamma}^T X_{\varepsilon, \gamma}
+\widehat{\Delta}_{\varepsilon, \gamma}+\varepsilon I)^{-1} -
(X_{\varepsilon, \gamma}^T X_{\varepsilon, \gamma} +\varepsilon
I)^{-1}\right) \right\|
  \\
& & = \left\|\left[(X_{\varepsilon, \gamma}^T X_{\varepsilon,
\gamma}+\varepsilon I) -\varepsilon I \right]
\left((X_{\varepsilon, \gamma}^T X_{\varepsilon, \gamma}
+\widehat{\Delta}_{\varepsilon, \gamma}+\varepsilon I)^{-1} -
(X_{\varepsilon, \gamma}^T X_{\varepsilon, \gamma} +\varepsilon
I)^{-1}\right) \right\|
\\
  &  & =  \left\| \left(I -\varepsilon (X_{\varepsilon, \gamma}^T
X_{\varepsilon, \gamma} +\varepsilon I)^{-1}\right) \left( \sum_{
i=1}^\infty (-1)^{i}
 \left(\widehat{\Delta}_{\varepsilon, \gamma}
(X_{\varepsilon, \gamma}^T X_{\varepsilon, \gamma} +\varepsilon
I)^{-1} \right)^{i}\right) \right\| \\
& &  \leq \left\| \left(I-\varepsilon (X_{\varepsilon, \gamma}^T
X_{\varepsilon, \gamma} +\varepsilon I)^{-1} \right) \right\|
\left\|\sum_{ i=1}^\infty (-1)^{i}
 \left(\widehat{\Delta}_{\varepsilon, \gamma}
(X_{\varepsilon, \gamma}^T X_{\varepsilon, \gamma} +\varepsilon
I)^{-1} \right)^{i} \right\| \\
& & \leq  \left\|\sum_{ i=1}^\infty (-1)^{i}
 \left(\widehat{\Delta}_{\varepsilon, \gamma}
(X_{\varepsilon, \gamma}^T X_{\varepsilon, \gamma} +\varepsilon
I)^{-1} \right)^{i} \right\|\\
& & \leq  \sum_{ i=1}^\infty \left\|
 \widehat{\Delta}_{\varepsilon, \gamma} \right\|^i  \left\| (X_{\varepsilon, \gamma}^T
X_{\varepsilon, \gamma} +\varepsilon I)^{-1} \right\|^i   \leq
\sum_{ i=1}^\infty \left(\|
 \widehat{\Delta}_{\varepsilon, \gamma} \| /\varepsilon
 \right)^i = \frac{\left(\|
 \widehat{\Delta}_{\varepsilon, \gamma} \| /\varepsilon
 \right)}{1- \left(\|
 \widehat{\Delta}_{\varepsilon, \gamma} \| /\varepsilon
 \right)}.
\end{eqnarray*}
Thus, by (\ref{limlim}), we have \begin{equation}\label{zerolim}
~~~~~~ \textrm{tr} \left( X_{\varepsilon, \gamma}^TX_{\varepsilon,
\gamma} \left[ \left( X_{\varepsilon, \gamma}^T X_{\varepsilon,
\gamma} +\widehat{\Delta}_{\varepsilon, \gamma}+\varepsilon I
\right)^{-1} - \left( X_{\varepsilon, \gamma}^T X_{\varepsilon,
\gamma} +\varepsilon I \right)^{-1} \right]\right)
  \to 0
\end{equation}
as $(\varepsilon, \gamma) \to 0$  and $\gamma /\varepsilon \to 0.$
By (\ref{P3}), (\ref{relation}) and (\ref{P1}),  we have
 \begin{eqnarray*} & & \textrm{rank}(X^*)  \geq  \phi_\varepsilon (X^*) \geq       \textrm{tr}(Y_{\varepsilon,
 \gamma})  =     \textrm{tr} \left(X_{\varepsilon, \gamma} \left(Z_{\varepsilon,
 \gamma}+\varepsilon I\right)^{-1} X_{\varepsilon, \gamma}^T \right)\\
 & & ~ =      \textrm{tr} \left(X_{\varepsilon, \gamma} \left(  X_{\varepsilon,
\gamma}^T X_{\varepsilon, \gamma} +\widehat{\Delta}_{\varepsilon,
\gamma}+\varepsilon I  \right)^{-1} X_{\varepsilon, \gamma}^T \right)\\
&  &  ~=     \textrm{tr} \left( X_{\varepsilon, \gamma} \left[
\left(X_{\varepsilon, \gamma}^T X_{\varepsilon, \gamma}
+\widehat{\Delta}_{\varepsilon, \gamma}+\varepsilon I  \right)^{-1}
- \left( X_{\varepsilon, \gamma}^T
X_{\varepsilon, \gamma}  +\varepsilon I \right)^{-1} \right] X_{\varepsilon, \gamma}^T \right) \\
& & ~~~~~~ +  \textrm{tr} \left(X_{\varepsilon, \gamma} \left(
X_{\varepsilon, \gamma}^T X_{\varepsilon, \gamma}  +\varepsilon I
\right)^{-1}
X_{\varepsilon, \gamma}^T\right) \\
& &  ~=  \textrm{tr} \left( X_{\varepsilon, \gamma}^T
X_{\varepsilon, \gamma} \left[ \left(X_{\varepsilon, \gamma}^T
X_{\varepsilon, \gamma} +\widehat{\Delta}_{\varepsilon,
\gamma}+\varepsilon I  \right)^{-1} - \left( X_{\varepsilon,
\gamma}^T X_{\varepsilon, \gamma} +\varepsilon I \right)^{-1}
\right] \right) + \phi_\varepsilon (X_{\varepsilon, \gamma}).
\end{eqnarray*}
which together with (\ref{zerolim}) and Lemma 2.5 implies that
\begin{eqnarray*}   & &    \textrm{rank(}X^*)  \geq     \limsup_{\varepsilon, \gamma\to 0,
\frac{\gamma}{\varepsilon} \to
 0} \textrm{tr}(Y_{\varepsilon, \gamma}) \geq  \liminf_{\varepsilon, \gamma\to 0,
\frac{\gamma}{\varepsilon} \to
 0} \textrm{tr}(Y_{\varepsilon, \gamma}) \\
   &   &  =      \liminf_{\varepsilon, \gamma\to 0,
\frac{\gamma}{\varepsilon} \to
 0} \{ \textrm{tr} ( X_{\varepsilon, \gamma}^T X_{\varepsilon, \gamma}
[ (X_{\varepsilon, \gamma}^T X_{\varepsilon, \gamma}
+\widehat{\Delta}_{\varepsilon, \gamma}+\varepsilon I  )^{-1}  - (
X_{\varepsilon, \gamma}^T X_{\varepsilon, \gamma} +\varepsilon I
)^{-1} ] )\\
& & ~~~~~~~~~~~~~~~   +  \phi_\varepsilon ( X_{\varepsilon, \gamma}) \} \\
  & &  =        \liminf_{\varepsilon,
\gamma\to 0, \frac{\gamma}{\varepsilon }\to
 0} \phi_\varepsilon ( X_{\varepsilon, \gamma})
 \geq   \Phi_0(X^0) =\textrm{rank}(X^0).
\end{eqnarray*}
Since $X^*$ is a minimum rank solution, all inequalities  above must
be equalities, and thus $X^0$ is a minimum rank solution, and
$\lim_{\varepsilon, \gamma\to 0, \frac{\gamma}{\varepsilon} \to
 0} \textrm{tr} (Y_{\varepsilon, \gamma}) =\textrm{rank}(X^*).$  ~ $  \endproof$

By Theorem 3.1, we may simply set $\gamma = \gamma (\varepsilon) $
as a function of $\varepsilon,$ for instance, $\gamma =\varepsilon^p
$ where $p>1 $ is a constant.    Then (\ref{phi-mini-C}) becomes the
problem below:
\begin{equation} \label{pppp} \begin{array}{cl} \textrm{Minimize}  & \textrm{tr}(Y) + \frac{1}{\gamma(\varepsilon)} \textrm{tr}(Z)\\
s.t. &   \left(
           \begin{array}{cc}
             Y & X \\
             X^T &  Z+ \varepsilon I \\
           \end{array}
         \right)
 \succeq  0,    \left(
                  \begin{array}{cc}
                    I & X  \\
                    X^T & Z \\
                  \end{array}
                \right) \succeq 0 ,
 ~~ X\in C,
\end{array}
\end{equation} which includes only one parameter.  An
immediate corollary from Theorem 3.1 is given as follows, which
shows that the minimum rank of (\ref{rank})
 can be obtained exactly by
solving (\ref{pppp}) with a suitable chosen parameter $\varepsilon.$

\textbf{Corollary 3.2.} \emph{Let $C \subset R^{m\times n} $  be a
bounded and  closed set, containing an element $X^*$ with the
minimum rank $r^*=rank(X^*)$ and the least F-norm.  Let $\gamma: (0,
\infty) \to (0, \infty) $ be a function satisfying
$\gamma(\varepsilon)/\varepsilon \to 0 $ as $\varepsilon \to 0.$  If
$(Y_\varepsilon, Z_\varepsilon, X_\varepsilon) $ is the optimal
solution of (\ref{pppp}), then  $tr(Y_\varepsilon) \leq r^*$ for all
$\varepsilon,$ and
$$ \lim_{\varepsilon\to 0}  \textrm{tr}(Y_\varepsilon)
 = r^* ,  ~~ \lim_{\varepsilon\to 0} \textrm{tr}(Z_\varepsilon)=
\|X^*\|_F^2,
$$ and the sequence $\{X_\varepsilon \}$ converges to the set of
   minimum rank solutions of (\ref{rank}). Moreover,  there exists a threshold $\delta>0$ such that
$ r^* = \left\lceil \textrm{tr}(Y_ \varepsilon)\right\rceil$  for
every  $\varepsilon \in (0, \delta].   $ }

From the above results, we see that a rank minimization problem can
be tractable under some conditions. We summarize this result as
follows.

\textbf{Corollary 3.3.} \emph{When $C$ is defined by linear
constraints (such as (\ref{CCC01}) and (\ref{CCC})), and if $C$
contains a minimum rank element with the least F-norm,    the rank
minimization problem (\ref{rank}) is equivalent to the SDP problem
(\ref{phi-mini-C}) by a suitable choice of the parameter $(\eta,
\varepsilon).$}

Note that the first term of the objective of (\ref{phi-mini-C}) is
to estimate $\textrm{rank}(X)$ and  the second term is to measure
the least F-norm.  So from Theorem 3.1 we may roughly say that under
some conditions minimizing $\textrm{rank}(X)$ over $C$ is equivalent
to minimizing $\textrm{rank} (X)+ \beta \|X\|_F^2$ over $C$  for
some $\beta. $ This is true,  as shown by the  next result below.

\textbf{Theorem 3.4.} \emph{Let the feasible set  be  of the form $
C= {\cal F} \cap \{X:  \gamma_1  \leq \|X\|_F \leq \gamma_2\}, $
where  $0< \gamma_1 \leq \gamma_2 $ are constants  and $ {\cal F}
\subseteq R^{m\times n}$  is a closed set.}

(i) \emph{The following two problems are equivalent in the sense
that they yield the same minimum rank solution:}
\begin{eqnarray} \label{fnorm-1}  &    \textrm{Minimize}  \left\{ \textrm{rank}
(X): ~ X\in C=  {\cal F} \cap\{X:  \gamma_1 \leq \|X\|_F \leq
\gamma_2\}\right\},
 & \\
& \label{fnorm-2}  ~~~~~~~~  \textrm{Minimize}  \left\{
\textrm{rank} (X) +  (1/\eta)  \|X\|_F^2: ~  X\in C=  {\cal F}
\cap\{X: \gamma_1 \leq \|X\|_F \leq \gamma_2\}\right\}, &
\end{eqnarray}
where $\eta > \gamma_2-\gamma_1 $ is a given number

 (ii) \emph{If ${\cal F} $ is a cone, then the
  set $C$ contains a minimum rank matrix $X^*$ with the least F-norm.}

\emph{Proof.}   (i) Assume that $X^*$ is a minimizer of
(\ref{fnorm-1}) with the minimum rank  $r^*,$ and assume that
$\widetilde{X }$ is an arbitrary minimizer of the problem
(\ref{fnorm-2}). We show that $\textrm{rank}(\widetilde{X})=r^*.$ In
fact, if this is not true, then $\textrm{rank}(\widetilde{X}) \geq
r^*+1 ,$ and thus
\begin{eqnarray}  \label {FFFF} \textrm{rank}(\widetilde{X})+
(1/\eta) \|\widetilde{X}\|_F^2 & \geq & \textrm{rank} (X^*)+1+
(1/\eta) \|\widetilde{X}\|_F^2 \nonumber\\  &=  &  \textrm{rank}
(X^*)+(1/\eta) \|X^*\|_F^2 + 1+ (1/\eta) \left(\|\widetilde{X}\|_F^2
-
  \|X^*\|_F^2\right)  \nonumber \\
  & > & \textrm{rank} (X^*)+ (1/\eta)
\|X^*\|_F^2,
\end{eqnarray}  where the last inequality above follows from the fact
$X^*, \widetilde{X}\in \{X: \gamma_1 \leq \|X\|_F\leq \gamma_2\}$
which implies that $ 1+ (1/\eta) \left(\|\widetilde{X}\|_F^2 -
  \|X^*\|_F^2\right) \geq 1+ (1/\eta)(\gamma_1-\gamma_2) >0
  $ by the choice of $\eta.$
Thus, (\ref{FFFF})  contradicts with the fact of $\widetilde{X}$
being a minimizer of (\ref{fnorm-2}).

(ii) Suppose that ${\cal F}$ is cone. Consider the F-norm
minimization problem:
$$   \textrm{Minimize}  ~\{ \|X\|_F^2:  X\in C= F\cap\{X:  \gamma_1  \leq \|X\|_F \leq \gamma_2\} \}.
$$
Since the feasible set of the  problem is  closed and  bounded, the
least F-norm solution, denoted by  $\overline{X},$ exists. Let $X^*$
be a minimum rank element in $C. $ Then $\|X^*\|_F\geq \|
\overline{X} \|_F \geq \gamma_1 >0 .$ Thus, there is a positive
number $1\geq \alpha>0$ such that $ \alpha \|X^*\|_F =
\|\overline{X}\|_F. $ Note that $\alpha X^* \in {\cal F} $ (since
${\cal F} $ is a cone), and that $\textrm{rank}(\alpha
X^*)=\textrm{rank}(X^*). $ Thus, $\alpha X^*$ is a  minimum rank
matrix with the least F-norm in $C.$  ~ $
\endproof$

Before we close this section, let us make some further comments on
the situation where $C$ is the intersection of a cone and a bounded
set defined by matrix norm, as discussed in Theorem 3.4.  This
situation does arise in the study of quadratic (in)equality systems
and quadratic optimization. First of all, it is worth pointing out
the following fact. Its proof is evident and omitted.

 \textbf{Theorem 3.5.}  \emph{Let ${\cal F}$  be a cone in $R^{m\times n},$
and let $ 0< \gamma_1\leq \gamma_2$ be two positive numbers. Then
the minimum rank $r^*$ of the rank minimization problem
\begin{equation} \label{cone}
  r^*= \min \left\{ rank(X): ~
   X\in C={\cal F}\cap \{X: \gamma_1\leq  \|X\| \leq
\gamma_2\} \right\}
\end{equation}
is independent of the choice of $\gamma_1, \gamma_2$ and the norm
$\|\cdot\|. $}

In another word,  no matter what matrix norms and the positive
numbers $\gamma_1, \gamma_2$ are used, the problem of the form
(\ref{cone}) yields the same minimum rank. So, in theory,  all these
rank minimization problems
  are equivalent.  From a computation point of
view, however, the choice of the norm $\|\cdot\|$ does matter. For
instance, when ${\cal F} $ is a subset of the positive semidefinite
cone, there are some benefits  of using the nuclear norm $\|X\|_*$
in (\ref{cone}). Since $\|X\|_* =\textrm{tr} (X) $ in positive
semidefinite cone, the constraint $ \gamma_1\leq \|X\|_*\leq
\gamma_2$ in this case coincides with the linear constraint
$\gamma_1\leq \textrm{tr}(X)\leq \gamma_2.  $ As a result, the
approximation counterpart, defined by (\ref{phi-mini-C}), of the
problem (\ref{cone}) is an SDP problem for this case, and hence it
can be solved efficiently. However, when the nuclear norm is used in
(\ref{cone}), the problem (\ref{cone}) may not satisfy the condition
of Theorem 3.1.

 When $C$ is defined by a
cone,  from Theorem 3.4 (ii) the problem (\ref{cone}) satisfies  the
condition of Theorem 3.1. However, when the F-norm is used, the
problem (\ref{phi-mini-C}) is not convex in general. To handle this
nonconvexity, we may consider the relaxation of (\ref{cone}). For
instance,  when ${\cal F}$ in (\ref{cone}) is a cone contained in
the positive semidefinite cone,  we define
  \begin{equation} \label{dddd}  \left\{\begin{array} {cc} \delta_1 =\min\{\textrm{tr}(X):    \gamma_1 \leq \|X\|_F \leq
\gamma_2, ~X \succeq 0 \}, \\
\delta_2 =\max\{\textrm{tr}(X): \gamma_1 \leq \|X\|_F \leq \gamma_2,
~ X \succeq 0 \} \end{array} \right. \end{equation} where $\gamma_1
>0.$   Clearly,
$\delta_1$ and $ \delta_2 $ exist and are positive. Thus the problem
(\ref{cone}) is relaxed to
$$
   l^* = \min \{ \textrm{rank}(X):  ~~ X\in C={\cal F}\cap \{X: \delta_1\leq  \textrm{tr}(X) \leq
\delta_2\}
 \}.
$$
When $F$ is defined by linear constraints, the approximation
counterpart (\ref{phi-mini-C}) of this relaxation problem is an SDP
problem.  Denote the optimal solution of this SDP problem by
$(Y_{\varepsilon, \eta}, Z_{\varepsilon, \eta}, X_{\varepsilon,
\eta}).  $ Then  by Theorem 3.1 it provides a lower bound for the
minimum rank of the above relaxation problem, and hence a lower
bound for the minimum rank of the original problem (\ref{cone}),
i.e.,  $ tr(Y_{\varepsilon, \eta}) \leq l^* \leq r^*.$

\section{Application to the system of quadratic equations}
Given a finite number of matrices $A_i\in S^n, i=1, ...m, $  we
consider the development of sufficient conditions for the following
assertion:
\begin{equation} \label{000} x^T A_i x=0, ~ i=1, ..., m
\Longrightarrow x=0,
\end{equation} i.e., $0$ is the only solution to (\ref{quadratic}).    At the first glance,
it seems that (1.5) and (4.1) have nothing to do with a rank
minimization problem.  In this section, however, we show that (4.1)
can be equivalently formulated as a rank minimization problem, based
on which we may derive some sufficient conditions for (\ref{000}) by
applying the approximation theory  developed in previous sections.
Note that system (\ref{quadratic}) can be written as $ \langle A_i,
xx^T\rangle=0, ~ i=1,\dots, m. $ Since $X=x x^T$ is either 0 (when
$x=0$) or a positive semidefinite rank-one matrix (when $x\not=0$),
it is natural to  consider the
 linear system:
\begin{equation}\label{rankone-M} \langle A_i, X\rangle=0, ~ i=1,\dots,
m, ~ X\succeq 0,
\end{equation}
which is a homogeneous system. The set $\{ X:  \langle A_i,
X\rangle=0, ~ i=1,\dots, m, ~X\succeq 0\}$ is a convex cone. It is
evident that \emph{the system (\ref{quadratic}) has a nonzero
solution if and only if the system (\ref{rankone-M}) has a rank-one
solution. In another word, $0$ is the only solution to
(\ref{quadratic}) if and only if (\ref{rankone-M}) has no rank-one
solution.} There are only two cases for the system (\ref{rankone-M})
with no  rank-one solution: either $X=0$ is the only matrix
satisfying (\ref{rankone-M}) or the minimum rank of the nonzero
matrices satisfying (\ref{rankone-M}) is greater than or equal to 2.
As a result, let us consider the following rank minimization
problem:
\begin{equation}\label{rank-optim} ~~~~~~~~ r^* = \min \left\{
\textrm{rank}(X): ~ \langle A_i, X\rangle=0, ~ i=1,\dots, m, ~
\delta_1 \leq \|X\|\leq \delta_2, ~X\succeq 0\right\},
\end{equation} where $0< \delta_1\leq \delta_2 $ are two given positive
constants. Clearly,  $X=0$ is the only matrix satisfying
(\ref{rankone-M}) if and only if the problem (\ref{rank-optim}) is
infeasible, in which case we set $r^*=\infty. $ It is also easy to
see that system (\ref{rankone-M}) has a  solution $X\not=0 $ if and
only if the problem (\ref{rank-optim}) is feasible, in which case
$r^*$ is finite and $ 1 \leq r^*\leq n.$  Thus for the problem
(\ref{rank-optim}),  we have either $r^* =\infty $ or $ 1\leq r^*
\leq n. $

From the above discussion, we immediately have the following result.

\textbf{Lemma 4.1.}   \emph{$0$ is the only solution to system
(\ref{quadratic})  if and only if  $r^*\geq 2 $ where $r^*$ is the
minimum rank of (\ref{rank-optim}).}

Thus developing a sufficient condition for (\ref{000})  can be
achieved by identifying the condition under which
   the minimum rank of (\ref{rank-optim}) is greater  than or equal
   to
   2. We follow this idea to establish some sufficient conditions for
(\ref{000}). By Theorem 3.5, the optimal value $r^*$ of
(\ref{rank-optim}) is independent of the choice of $\delta_1,
\delta_2$ and $\|\cdot\|. $ Thus Lemma 4.1 holds for any given
$0<\delta_1\leq \delta_2$ and any prescribed matrix norm in
(\ref{rank-optim}). So  we have a freedom to choose $\delta_1,
\delta_2 $ and the matrix norm in (\ref{rank-optim}) without
affecting the  value of $r^* $  in (\ref{rank-optim}). Thus, by
setting $\delta_1=\delta_2=1$ for simplicity and using the $F$-norm
in (\ref{rank-optim}), we have the problem
\begin{equation}\label{rank-optim-F}  ~~~~~~ r^* = \min \left\{
\textrm{rank}(X):  ~ \langle A_i, X\rangle=0, ~ i=1,\dots, m, ~
\|X\|_F= 1, ~X\succeq 0\right\}.
\end{equation}
By Theorem 3.4(ii),  the feasible set of this problem contains a
minimum rank solution with the least F-norm (which is equal to 1 for
this case). From Theorem 3.1 and its corollary, the rank
minimization (\ref{rank-optim-F}) can be approximated by the
following continuous optimization problem (as $(\eta, \varepsilon)
\to 0 $ and $\eta/\varepsilon \to 0$):
\begin{eqnarray} \label{PPPP}   & \textrm{Minimize }  &  \textrm{tr}(Y) +  (1/\eta) \textrm{tr}(Z)  \nonumber\\
 & s.t. &    \left(
           \begin{array}{cc}
             Y & X \\
             X &  Z+ \varepsilon I \\
           \end{array}
         \right)
 \succeq  0,    \left(
                  \begin{array}{cc}
                    I & X  \\
                    X & Z \\
                  \end{array}
                \right) \succeq 0,\\
  &  &    \langle A_i, X \rangle = 0, ~ i=1, ..., m,  ~   \|X\|_F=1, ~ ~ X\succeq
  0. \nonumber
\end{eqnarray}
(All results later in this section can be stated without involving
the parameter $\eta$ by setting,  for instance,  $\eta=
\varepsilon^2$ for the simplicity). By Corollary 3.2, the first term
of the objective in the above problem provides a lower bound for the
minimum rank of (\ref{rank-optim-F}).   However, the constraint
$\|X\|_F=1 $ makes the problem (\ref{PPPP}) difficult to be solved
directly.  So let us consider a relaxation  of this constraint.
Similar to (\ref{dddd}), we define two constants:
\begin{equation}\label{constant} ~~~~~ \delta_1= \min \{\textrm{tr}(X): \|X\|_F =1,
X\succeq 0\}, ~~ \delta_2= \max \{\textrm{tr}(X): \|X\|_F =1,
X\succeq 0\}.
\end{equation}
It is easy to verify that $\delta_1=1 $ and $\delta_2=\sqrt{n}.$ In
fact, in terms of eigenvalues of $X$, the above two extreme problems
are nothing but minimizing and maximizing, respectively, the
function $\sum_{i=1}^n \lambda_i$ subject to $ \sum_{i=1}^n
\lambda_i^2 =1, \lambda_i\geq 0, i=1,...,n . $ The optimal values of
these two problems are $1$ and $\sqrt{n},$ respectively. Therefore,
we conclude that
$$\{X: ~\|X\|_F=1, ~X\succeq 0\} \subseteq \{X: ~ 1 \leq \textrm{tr}(X)\leq \sqrt{n},
~X\succeq 0\}.$$ Thus, the following SDP problem is a relaxation of
(\ref{PPPP}):
\begin{equation} \label{RRRR} \begin{array}{cl} \textrm{Minimize}  & \textrm{tr}(Y) + (1/\eta)  \textrm{tr}(Z)\\
& \\ s.t. &   \left(
           \begin{array}{cc}
             Y & X \\
             X &  Z+ \varepsilon I \\
           \end{array}
         \right)
 \succeq  0,    \left(
                  \begin{array}{cc}
                    I & X  \\
                    X & Z \\
                  \end{array}
                \right) \succeq 0,\\
 & \\
  &     \langle A_i, X \rangle = 0, ~ i=1, ..., m,  ~ 1\leq \textrm{tr}(X)\leq \sqrt{n}, ~ ~ X\succeq
  0.
\end{array}
\end{equation}
The optimal value of (\ref{RRRR}) is a lower bound for that of
(\ref{PPPP}). We now derive out the dual problem of
  (\ref{RRRR}), which will be used to develop the sufficient
  condition
  for (\ref{000}). We will make  use of the following lemma.

\textbf{Lemma 4.2.} \emph{If the SDP problem is of the form
 \begin{equation} \label{Primal} \min \{\langle
C_0, W \rangle : \langle C_i, W\rangle = b_i, i=1, ..., l, ~
\delta_1\leq \langle C, W \rangle \leq \delta_2, ~ W \succeq 0, \}
\end{equation}  then its dual problem is given by  \begin{equation}
\label{dual} \max \left\{ b^T y + \delta_1 t_1 +\delta_2 t_2:   ~
\sum_{i=1}^l y_i C_i +(t_1 +t_2) C \preceq C_0, ~t_1\geq 0, ~t_2\leq
0 \right\}, \end{equation} where $b=(b_1, ..., b_l)^T.$}

\emph{ Proof.}  For the  standard  SDP problem $\min\{\langle P,
W\rangle: ~ \langle P_i, W\rangle = b_i, ~i=1, ..., l, ~ W\succeq
0\},$ it is well known that its dual problem is given by $\max \{
b^T y:  \sum_{i=1}^l y_i P_i \preceq P\}.$ By
 transforming (\ref{Primal}) into  the standard form and applying this fact, it is
easy to verity that the dual problem of (\ref{Primal}) is indeed
given by (\ref{dual}). The detail is omitted. ~ $\endproof $

To obtain the dual problem of (\ref{RRRR}), let us rewrite the
problem (\ref{RRRR}) as the form of (\ref{Primal}).  Notice that the
 positive semidefinite conditions in
   (\ref{RRRR}) are equivalent to \begin{equation} \label{WWWW} \left(\begin{array}{ccccc}
                                     X  &  0  &   0  & 0 & 0\\
                                     0  &  I  &  X   & 0 & 0  \\
                                     0  &  X  &  Z   & 0 & 0   \\
                                     0  &  0  &  0   & Y & X  \\
                                     0  &  0  &  0   & X & Z+\varepsilon I
                                     \end{array}
 \right)  \succeq 0 . \end{equation}
 Let $ E^{(k, l)} \in S^ {5n\times 5n} $ ($k,l=1,..., 5n)$ denote
the symmetric matrices with $(k, l)$th entry = $(l, k)$th entry $=
1$ and zero elsewhere. When $k=l$, $ E^{(k, k)} $ denotes the matrix
with ($k,k)$th entry 1 and all other elements 0. Clearly, we have
$E^{(l,k)}=E^{(k,l)}$ for any $(k,l).$  Note that  for any matrix
$W=(w_{i,j}) \in S^{5n\times 5n},$  it can be represented as
$W=\sum_{k=1}^{5n} \sum_{l=k}^{5n} w_{k,l}E^{(k,l)},$ and   $
\langle E^{(k, l)}, W \rangle  = w_{k,l}+w_{l,k} = 2w_{k,l}$ for
$k\not=l,$  and  $ \langle E^{(k, k)}, W \rangle  = w_{k,k}. $ In
terms of $E^{(\cdot, \cdot)},$ the condition (\ref{WWWW})
 can be written as the
following set of constraints:
\begin{eqnarray} W & \succeq  & 0 ,  \nonumber \\
 \langle E^{(i, ~n+j)}, W \rangle &  = & 0,  ~~i = 1, ..., n, j=1,...,
4n, \label{F1} \\
 \langle E^{(n+i, ~3n+j)}, W \rangle & = & 0,  ~~i,j = 1, ..., 2n,
\label{F2} \\
  \langle E^{(n+i, ~n+j)}, W \rangle & = &  0, ~ i=1, ..., n-1, j=i+1,...,
  n,
  \label{F3} \\
  \langle E^{(n+i, ~n+i)}, W \rangle & = &  1, ~ i=1,..., n,
  \label{F3b}\\
  \langle E^{(i,j)}-E^{(n+i, ~2n+j)}, W \rangle & = & 0, ~~ i=1,..., n-1, ~ j=i+1, ...,
  n,\label{F4} \\
   \langle E^{(j,i)}-E^{(n+j, ~2n+i)}, W \rangle  & = & 0, ~i=1,...,
n-1, ~j=i+1, ..., n, \label{F4b}\\
  \langle 2 E^{(i,i)}-E^{(n+i, ~2n+i)}, W \rangle  & = & 0, ~~ i =1, ...,
  n,\label{F5} \\
  \langle E^{(n+i,2n+j)}-E^{(3n+i, ~4n+j)}, W \rangle & = & 0, ~~ i,j=1, ..., n, \label{F6} \\
  \langle E^{(4n+i, ~4n+j)}- E^{(2n+i, ~2n+j)} , W \rangle  & = & 0,  ~~i=1,..., n-1,   j =i+1,..., n
  ,\label{F7} \\
  \langle E^{(4n+i, ~4n+i)}- E^{(2n+i, ~2n+i)} , W \rangle &  = & \varepsilon,  ~~i=1, ..., n ,
  \label{F8}
  \end{eqnarray}
  where (\ref{F1})
and (\ref{F2}) represent the zero blocks in the matrix of
(\ref{WWWW}), conditions (\ref{F3}) and (\ref{F3b}) describe the
block `$I$' (the identity matrix), conditions (\ref{F4})-(\ref{F6})
represent the `$X$' blocks, and (\ref{F7}) and (\ref{F8}) describe
the relation between the blocks `$Z$' and `$Z+\varepsilon I$'
therein.   In terms of $W\in S^{5n\times 5n},$  the
                                   equality $\langle A_i, X\rangle =0$ in (\ref{RRRR}) can be written as
$\langle P_i, W\rangle =0,  $  the inequality $  1\leq
\textrm{tr}(X)\leq \sqrt{n} $   can be represented  as $ 1 \leq
\langle  P_0 , W \rangle \leq \sqrt{n}, $ and the objective of
(\ref{RRRR}) can be written as $\langle P, W\rangle $ where  $P_i,
P_0, P\in S^{5n\times 5n}$ are given by $$ P_i = \left[
\begin{array}{cc}
                                     A_i &  0\\
                                     0&  0 \\
                                   \end{array}
 \right],  ~ P_0 =
\left(\begin{array}{cc}
                                     I     &  0 \\
                                     0  &  0
                                   \end{array} \right), ~ P=
\left(\begin{array}{ccc}
                                     0 &    &      \\
                                       &  \left(\begin{array} {cc} 0 &     \\
                                          & \frac{1}{\eta} I  \end{array} \right) &  \\
                                       &    &   \left(\begin{array}{cc}   I  &  \\
                                              &  0
                                   \end{array}  \right) \end{array}
                                   \right),
                                   $$
where $I$ is the $n\times n $ identity matrix. Thus,  (\ref{RRRR})
can be written as  the following  SDP problem:
\begin{eqnarray}  \label {SDP-02} & \min  & \left\langle
P,
W \right\rangle  \nonumber \\
& \textrm{s.t.} & \langle E^{(i, ~n+j)}, W \rangle =0,  ~~i = 1,
..., n, j=1,...,
4n,  \nonumber \\
& & \langle E^{(n+i, ~3n+j)}, W \rangle =0,  ~~i,j = 1, ..., 2n,  \nonumber \\
& &  \langle E^{(n+i, ~n+j)}, W \rangle =  0,   ~ i=1, ..., n-1, ~j=i+1, ..., n, \nonumber \\
& &  \langle E^{(n+i, ~n+i)}, W \rangle =  1, ~  i=1, ..., n,  \nonumber \\
& &  \langle E^{(i,j)}-E^{(n+i, ~2n+j)}, W \rangle =0, ~i=1,...,
n-1, ~j=i+1,
..., n,    \\
& &  \langle E^{(j,i)}-E^{(n+j, ~2n+i)}, W \rangle =0, ~i=1,...,
n-1, ~j=i+1,
..., n,   \nonumber \\
& & \langle 2 E^{(i,i)}-E^{(n+i, ~2n+i)}, W \rangle =0, ~~ i =1,
..., n,  \nonumber \\
 & & \langle E^{(n+i,2n+j)}-E^{(3n+i, ~4n+j)}, W \rangle =0, ~~ i,j=1, ..., n,
 \nonumber \\
& &  \langle E^{(4n+i, ~4n+j)}- E^{(2n+i, ~2n+j)}, W \rangle =0,
~~i=1,..., n-1,~ j =i+1, ..., n ,  \nonumber \\
& &  \langle E^{(4n+i, ~4n+i)}- E^{(2n+i, ~2n+i)}, W \rangle
=\varepsilon,  ~~i=1, ..., n ,  \nonumber \\
& & \left\langle P_i,
W \right\rangle  =  0, ~~ i=1, ..., m, \nonumber \\
& &  1 \leq \langle  P_0 , W \rangle \leq  \sqrt{n}, \nonumber
\\
 & &  W \succeq 0 \nonumber
\end{eqnarray}
 which  is of the
form (\ref{Primal}). By Lemma 4.2, its dual problem  is given by
\begin{eqnarray*}  \label {DDDI} & \max &
\sum_{i=1}^n \alpha_i+ \sum_{i=1}^n \varepsilon \beta_i +   t_1+ \sqrt{n} t_2 \nonumber \\
& \textrm{s.t.}& \nonumber \\
 & &
 \sum_{i=1}^n\sum_{j=1}^{4n} \rho_{ij} E^{(i, n+j)} +  \sum_{i,j=1}^{2n} \rho'_{ij} E^{(n+i, 3n+j)}
 + \sum_{i=1}^{n-1}\sum_{ j=i+1 }^n     \rho''_{ij} E^{(n+i, n+j)}
  \nonumber \\
   & & + \sum_{i=1}^n    \alpha_i  E^{(n+i, n+i)} + \sum_{i=1}^{n-1}\sum_{ j=i+1 }^n \left[ \xi_{ij} (E^{(i,j)}-E^{(n+i, 2n+j)})
 +   \xi'_{ij} (E^{(j,i)}-E^{(n+j, 2n+i)}) \right]   \\
 & &   + \sum_{i=1}^{n} \eta_{i}
 (2E^{(i,i)}-E^{(n+i, 2n+i)})
  + \sum_{i,j=1 }^{n} \theta_{ij}  (E^{(n+i,2n+j)}-E^{(3n+i,
 4n+j)}) \nonumber\\
  & &  + \sum_{ i=1}^{n-1}\sum_{ j=i+1 }^{n} \theta'_{ij}  (E^{(4n+i,
 4n+j)}- E^{(2n+i,2n+j)}) + \sum_{i=1 }^{n} \beta_i (E^{(4n+i,
 4n+i)}- E^{(2n+i,2n+i)}) \nonumber\\
 & &   +  \sum_{i=1}^m y_i P_i +t_1 P_0+t_2 P_0 \preceq P, \nonumber\\
 & & t_1\geq 0, ~ t_2\leq 0. \nonumber
\end{eqnarray*}
  By the structure of $ P, P_0 , E^{(\cdot,\cdot)}$'s
and  $P_i$'s, the above problem can be written as {\small
\begin{eqnarray} \label {DDDD} & ~~~\max   &
\sum_{i=1}^n \alpha_i+ \sum_{i=1}^n \varepsilon \beta_i +
t_1+\sqrt{n}  t_2 ~~\left( = \textrm{tr}(\Phi) -\varepsilon
\textrm{tr}(Q) + t_1+\sqrt{n} t_2\right)
   \nonumber \\
&\textrm{ s.t.}& \nonumber \\
 & &
  \left(
    \begin{array}{ccccc}
      V+V^T +\sum_{i=1}^m y_i  A_i+(t_1+t_2)I & U_1 & U_2 & U_3 & U_4 \\
      U_1^T & \Phi & \Theta-V & U_5 & U_6 \\
      U_2^T & \Theta^T-V^T & Q-\frac{1}{\eta}I & U_7 & U_8 \\
      U_3^T & U_5^T & U_7^T &  -I & -\Theta \\
      U_4^T & U_6^T & U_8^T & -\Theta^T & -Q \\
    \end{array}
  \right) \preceq 0, \\
 & & t_1\geq 0, ~ t_2\leq 0, \nonumber
\end{eqnarray} }
where $\alpha_i$'s and $-\beta_i$'s  are the diagonal entries of
 $\Phi$ and $Q$, respectively, i.e.,
 $ \textrm{diag} (\Phi)= (\alpha_1, ..., \alpha_n),$  and $ \textrm{diag} (Q)= (-\beta_1,
..., -\beta_n).$   Thus the objective of the above problem can be
written as $\textrm{tr}(\Phi) -\varepsilon \textrm{tr}(Q) +
t_1+\sqrt{n}t_2. $ All blocks in the above matrix are $n\times n $
submatrices. Also, note that ({\ref{DDDD}) is always feasible and
satisfies the Slater's condition, for instance, $(\Theta=V=0,
\Phi=-I, Q =\frac{1}{2\eta} I, t_1=1, t_2=-2, y_i=0 $ for all $
i=1, ..., m, $ and $ U_i=0 $ for all $ i=1, ..., 8) $ is a strictly
feasible point.

\textbf{Theorem 4.3. } \emph{  If there exist $(\eta,
\varepsilon)>0$ and $t_1, t_2 , \mu_i, i=1,...,m $ and matrices
$\Phi, Q \in S^{n\times n}, V, \Theta \in R^{n\times n} $ and
$M_i\in R^{n\times n }, i=1,...,8$ such that the following
conditions hold {\small
\begin{eqnarray}\label{MM}  & \left\lceil \textrm{tr}(\Phi) -\varepsilon
\textrm{tr}(Q) +  t_1+ \sqrt{n} t_2
-\frac{1}{\eta} \right\rceil \geq 2, ~~ t_1\geq 0, ~ t_2\leq 0, & \\
&  \label{MMM} ~~~~~~~\left(
    \begin{array}{ccccc}
     \sum_{i=1}^m \mu_i  A_i - (t_1+t_2)I - (V+V^T)  & M_1 & M_2 & M_3 & M_4 \\
      M_1^T & -\Phi & V- \Theta & M_5 & M_6 \\
      M_2^T & V^T-\Theta ^T &  \frac{1}{\eta}I -Q & M_7 & M_8 \\
      M_3^T & M_5^T & M_7^T &  I  & \Theta \\
      M_4^T & M_6^T & M_8^T & \Theta^T & Q \\
    \end{array}
  \right)\succeq 0, & \end{eqnarray} }
then $0$ is the only solution to the quadratic equation
(\ref{quadratic}).}

\emph{Proof.}   Let $X^*$ be the minimum rank solution of
(\ref{rank-optim-F}) with the least norm $\|X^*\|_F=1.$ Let
$(Y_{\eta,\varepsilon}, X_{\eta, \varepsilon}, Z_{\eta,
\varepsilon}) $ be the optimal solution to (\ref{PPPP}), by Theorem
3.1, we have $r^* \geq \lceil Y_{\eta,\varepsilon} \rceil$ for every
$(\eta, \varepsilon)>0,$ where $r^*$ is the minimum rank of
(\ref{rank-optim-F}).  Since (\ref{RRRR}) is a relaxation of
(\ref{PPPP}), the optimal value of (\ref{RRRR}), denoted by
$v^*(\eta, \varepsilon) $, provides a lower bound for that of
(\ref{PPPP}), i.e.,
\begin{equation} \label{SC1}
\textrm{tr}(Y_{\eta, \varepsilon})+(1/\eta)\textrm{tr}(Z_{\eta,
\varepsilon}) \geq v^*(\eta, \varepsilon) ,
\end{equation} which holds for any given $(\eta, \varepsilon)>0. $
Note that  (\ref{DDDI}) is the dual problem of (\ref{RRRR}). If the
conditions (\ref{MM}) and (\ref{MMM}) hold, then for this $(\eta,
\varepsilon),$ the point $(t_1, t_2, y_i=-\mu, i=1, ..., m, \Phi, V,
\Theta, U_i= -M_j, j=1, ..., 8)$ is feasible to the dual problem
(\ref{DDDI}). Thus, by duality theory we have
 \begin{equation} \label{SC2}  v^*(\eta, \varepsilon) \geq  \textrm{tr}(\Phi) -\varepsilon \textrm{tr}(Q) +  t_1+\sqrt{n} t_2
.\end{equation} Notice that  $(Y^*,Z^*,X^*),$ where $Y^* =
X^*((X^*)^TX^*+\varepsilon I)^{-1} (X^*)^T$ and $ Z^* = (X^*)^TX^*,
$ is a feasible point of (\ref{PPPP}). Thus \begin{equation}
\label{SC3}  \textrm{tr}(Y_{\eta, \varepsilon})+(1/\eta)
\textrm{tr}(Z_{\eta,\varepsilon}) \leq \textrm{tr}(Y^*)+ (1/\eta)
\textrm{tr}(Z^*) = \phi_\varepsilon (X^*)+ (1/\eta),
\end{equation} where the last equality follows from that $ \textrm{tr}(Y^*) =
\phi_\varepsilon (X^*)$ and $\textrm{tr}(Z^*) = \|X^*\|_F^2=1.$
 Combining
(\ref{SC1}), (\ref{SC2}) and (\ref{SC3}) yields $$ \phi_\varepsilon
(X^*)+ (1/\eta) \geq \textrm{tr}(\Phi)-\varepsilon \textrm{tr}(Q)+
t_1+\sqrt{n} t_2 . $$   This together with (\ref{relation}) implies
that $ \textrm{rank}(X^*) \geq \textrm{tr}(\Phi)- \varepsilon
\textrm{tr}(Q)+ t_1+\sqrt{n} t_2 -(1/\eta).$ Thus, under the
conditions (\ref{MM}) and (\ref{MMM}), we see that
$$ r^*= \textrm{rank}(X^*) \geq \left\lceil \textrm{tr}(\Phi)-\varepsilon \textrm{tr}(Q)+  t_1+\sqrt{n} t_2
-(1/\eta)\right\rceil \geq 2. $$ By Lemma 4.1, we conclude that
(\ref{000}) holds, i.e., $0$ is the only solution to
(\ref{quadratic}). ~ $ \endproof. $

 From the above result,  a number
of sufficient conditions stronger than (\ref{MM})-(\ref{MMM}) can be
obtained. For example, we have the following corollary.

 \textbf{Corollary 4.4.} \emph{Let $A_i\in S^n, i=1, ...,m$ be a given set
 of
matrices. If there exist $(\eta,\varepsilon) >0, t_1, t_2, \mu_1,
..., \mu_m \in R, $ $Q , \Phi \in S^{n\times n}$ and $V, \Theta\in
R^{n\times n}$  such that
\begin{eqnarray} \left\lceil tr(\Phi) - \varepsilon tr(Q) +
t_1+\sqrt{n}t_2 - 1/\eta \right\rceil \geq 2, ~~~t_1 & \geq & 0,
~t_2\leq 0, \label{SS1} \\
\left(
    \begin{array}{cc}
      -\Phi & V- \Theta \\
      V^T-\Theta^T  &  \frac{1}{\eta}I -Q
      \end{array} \right) \succeq 0, ~~
     \left(\begin{array}{cc}
           I  & \Theta \\
           \Theta^T & Q
    \end{array}
  \right)  & \succeq & 0,  \label{SS2} \\
  \sum_{i=1}^m \mu_i  A_i - (t_1+t_2)I - (Y+Y^T)
 & \succeq  & 0, \label{SS3}
  \end{eqnarray}
     then $0$ is the only solution to the  system (\ref{quadratic}). }

We  now point out that   (\ref{positive}) implies
(\ref{SS1})-(\ref{SS3}).  Let $\eta >0$ be a given number.
  If $\sum_{i=1}^m t_i A_i \succ 0 $ for some $t_i,
 i=1,...,m,$ then we choose $\mu_i= \alpha t_i$ where $\alpha $
 can be any large positive number such that $\sum_{i=1}^m \mu_i A_i \succ
 t_1 I
 $ where $t_1=2+\frac{1}{\eta}.$ Then
 conditions (\ref{SS1})-(\ref{SS3}) hold with $V= \Phi= Q =\Theta=0$ and $t_2=0. $   Thus, the known condition
 (\ref{positive}) indeed implies (\ref{SS1})-(\ref{SS2}).  For $m=2 $ and $n\geq 3$, since the condition (\ref{000})
  is equivalent to $\mu_1 A_1+\mu_2 A_2 \succ 0,$ the
   sufficient conditions in Theorem 4.4 and Corollary 4.5 are also
   necessary conditions for (\ref{000}).
 When $n=2$ or $m\geq 3$,
  if  an example can be given to show that our
  sufficient conditions do not imply (\ref{positive}), then the open
  `Problem 13' in \cite{HU07} would be addressed.  We conjecture that
  our sufficient conditions are indeed mild ones for (\ref{000}).

\textbf{Remark 4.5.} To get more simple sufficient conditions for
(\ref{000}), we may continue to reduce the freedom of the variables
in  (\ref{SS1})-(\ref{SS3}). For instance, (\ref{SS2}) can be
replaced by a stronger version like $ \Phi \preceq 0,
\frac{1}{\varepsilon^2} I \succeq Q, Q \succeq Y^T Y $
  without involving  the matrix $\Theta.$
It is also worth stressing that  checking the new sufficient
  conditions developed in this section can be achieved by solving an SDP problem. For instance,
 if the optimal value of the SDP problem
  (\ref{DDDI}) is greater than $\frac{1}{\eta} +1,$ then the conditions (\ref{MM})-(\ref{MMM})
  hold. Similarly,
  if the optimal value of the SDP problem (\ref{DDDI}) with $M_i=0, i=1,
  ...,8$ is greater than $\frac{1}{\eta} +1, $  then
  the conditions (\ref{SS1})-(\ref{SS3}) hold.

\section{Conclusions}
Since $\textrm{rank}(X)$ is a discontinuous function with an integer
value,  this makes the rank minimization problem hard to be solved
directly. In this paper, we have presented a generic approximation
approach for rank minimization problems through  the approximation
function $\phi_\varepsilon(X). $ In particular,  we have shown that
when the feasible set is bounded the rank minimization problem can
be approximated to any level of accuracy by a nonlinear SDP problem
or a linear bilevel SDP problem with a special structure. To obtain a
tractable approximation of the rank minimization with linear
constraints,  the approximation model (\ref{phi-mini-C}) is
introduced, and is proved to be efficient for locating the minimum
rank solution of the problem if the feasible set  contains a minimum
rank element with the least F-norm. In this case, the rank
minimization problem is   equivalent to an SDP problem. This theory
was applied to a system of quadratic equations which can be
formulated as a rank minimization. Based on its approximation
counterpart,  we have developed some sufficient conditions for such
a  system with only zero solution.

\end{document}